\documentclass[12pt]{article}

\usepackage{amssymb,amsfonts,amsthm,amsmath}

\makeatletter \@addtoreset{equation}{section}
\@addtoreset{theequation}{thesection} \makeatother

\newtheorem{theorem}{Theorem}[section]
\newtheorem{lemma}{Lemma}[section]

\theoremstyle{remark}
\newtheorem{remark}{Remark}[section]
\newtheorem{example}{Example}[section]

\DeclareMathOperator{\ran}{Ran} \DeclareMathOperator{\Ker}{Ker}
\DeclareMathOperator{\Id}{Id} \DeclareMathOperator{\dom}{Dom}

\begin{document} 
 
 \begin{center}{\Large \bf 
Image of the Spectral Measure of a Jacobi Field\\ and the Corresponding Operators
 \\[5mm]\large  Yurij M. Berezansky, Eugene W. 
Lytvynov\\ and Artem~D.~Pulemyotov}\end{center}













\begin{abstract}
By definition, a Jacobi field $J=(J(\phi))_{\phi\in H_+}$
is a family of commuting selfadjoint three-diagonal operators in
the Fock space $\mathcal F(H)$. The operators $J(\phi)$ are
indexed by the vectors of a real Hilbert space $H_+$. The spectral
measure $\rho$ of the field $J$ is defined on the space $H_-$ of
functionals over $H_+$. The image of the measure $\rho$ under a
mapping $K^+:T_-\to H_-$ is a probability measure $\rho_K$ on
$T_-$. We obtain a family $J_K$ of operators whose spectral
measure is equal to $\rho_K$. We also obtain the chaotic
decomposition for the space $L^2(T_-,d\rho_K)$.
\end{abstract}

\vspace{2mm}

\vspace{3mm}


\section{Introduction}\label{introduction}

Consider a real Hilbert space $H$ and the corresponding symmetric
Fock space \begin{align}\label{intro Fock}\mathcal
F(H)=\bigoplus_{n=0}^\infty \mathcal F_n(H).\end{align} Let
\begin{align*}H_-\supset H\supset H_+\end{align*} be a
rigging of $H$ with the quasinuclear embedding $H_+\hookrightarrow
H$. Consider a Jacobi field $J=(\tilde J(\phi))_{\phi\in H_+}$ in
the space $\mathcal F(H)$. By definition, a Jacobi field is a
family of commuting selfadjoint operators which have a
three-diagonal structure with respect to the
decomposition~(\ref{intro Fock}). These operators are assumed to
linearly and continuously depend on the indexing parameter
$\phi\in H_+$. The concept of a Jacobi field was studied in \cite{YuBVLEL95}, \cite{EL95}, \cite{YuB98}, \cite{YuB98IEO},
\cite{YuB98MF}, and~\cite{YuB00}.

The above-mentioned papers provide the expansion of the Jacobi
field $J$ in generalized joint eigenvectors. The corresponding
Fourier transform appears to be a unitary operator between the
Fock space $\mathcal F(H)$ and the space $L^2(H_-,d\rho)$. The
measure $\rho$ on $H_-$ is called the spectral measure of $J$.
Note that the Jacobi field with the Gaussian spectral measure is
the classical free field in quantum field theory. The Jacobi field
with the Poisson spectral measure was actually discovered in
\cite{RHKP84} and~\cite{DS84}.

Jacobi fields are actively utilized in non-Gaussian white noise
calculus and the theory of stochastic processes,
see~\cite{YuBVLEL95}, \cite{EL95}, \cite{YuB98IEO}, \cite{YuB00},
\cite{YuB00IDA}, \cite{YuKEL00}, \cite{YuC02}, \cite{YuBELDM03},
\cite{YuBDM03}, \cite{EL02}, \cite{ELsubm}, \cite{DM03},
\cite{EL04}, and also~\cite{DNWS00} and~\cite{WS00}. Other
applications are to the integration of nonlinear
difference-differential equations, see~\cite{YuB98}. In the case
of a finite-dimensional $H$, the theory of Jacobi fields is
closely related to some results in~\cite{MGAK91}, \cite{MGAK94},
and~\cite{CDYX01}.

The problem of finding an operator family with a given spectral
measure often arises in applications. In some situations, the
given measure is equal to the image of the spectral measure of a
Jacobi field under a certain mapping. More precisely, let $\rho$
be the spectral measure of the field $J$. Consider a mapping
$K^+:H_-\to T_-$ with $T_-$ being a certain Hilbert space. This
mapping takes $\rho$ to the measure $\rho_K$ on $T_-$. Our paper
aims to find a family $J_K$ of operators whose spectral measure
equals  $\rho_K$. In other words, we track the changes of the
Jacobi field caused by mapping its spectral measure. Noteworthily,
if $K^+$ is an invertible operator, then $J_K$ appears to be
isomorphic to the initial family $J$.

We also study the chaotic decomposition of the space
$L^2(T_-,d\rho_K)$, which is derived through the orthogonalization
of polynomials on $T_-$.

Throughout this paper, we assume $K^+$ to be a bounded operator
with \linebreak $\Ker(K^+)=\{0\}$. We will also assume $\ran(K^+)$
to be dense in $T_-$. This assumption is not essential because the
measure $\rho_K$ is lumped on $\ran(K^+)$, and we can always
replace $T_-$ with the closure of $\ran(K^+)$ in $T_-$.

\section{Preliminaries}\label{sec_DPJF}

Let $H$ be a real separable Hilbert space. The corresponding
symmetric Fock space is defined as
\begin{align*}\mathcal F(H)=\bigoplus_{n=0}^\infty \mathcal
F_n(H)\end{align*} and consists of sequences
$\Phi=(\Phi_n)_{n=0}^\infty$, $\Phi_n\in\mathcal
F_n(H)=H_c^{\hat{\otimes}n}$, ($H_c$ being the complexification of
$H$ and $\hat\otimes$  denoting symmetric tensor product). The
finite vectors $\Phi=(\Phi_1,\ldots,\Phi_n,0,0,\ldots)\in\mathcal
F(H)$ form a linear topological space $\mathcal
F_{fin}(H)\subset\mathcal F(H)$. The convergence in $\mathcal
F_{fin}(H)$ is equivalent to the uniform finiteness and
coordinatewise convergence. The vector
$\Omega=(1,0,0,\ldots)\in\mathcal F_{fin}(H)$ is called vacuum.

Let \begin{align}\label{rigging_H}H_-\supset H\supset
H_+\end{align} be a rigging of $H$ with real Hilbert spaces $H_+$
and $H_-=(H_+)'$ (hereafter, $X'$ denotes the dual of the space
$X$). We suppose the inequality $\|\cdot\|_{H_+}\ge\|\cdot\|_H$ to
hold for the norms. We also suppose the embedding
$H_+\hookrightarrow H$ to be quasinuclear. The pairing
in~(\ref{rigging_H}) can be extended naturally to a pairing
between $\mathcal F_n(H_+)$ and $\mathcal F_n(H_-)$. The latter
can, in turn, be extended to a pairing between $\mathcal
F_{fin}(H_+)$ and $(\mathcal F_{fin}(H_+))'$. In what follows, we
use the notation $\langle\cdot,\cdot\rangle_H$ for all of these
pairings. Note that $(\mathcal F_{fin}(H_+))'$ coincides with the
direct product of the spaces $\mathcal F_n(H_-)$, $n\in\mathbb
Z_+$.

\subsection{Definition of a Jacobi field}\label{defin JF}

In the Fock space $\mathcal F(H)$, consider the family $(\mathcal
J(\phi))_{\phi\in H_+}$ of operator-valued Jacobi matrices
\begin{align*}
\mathcal J(\phi)=\left(
\begin{array}{cccccc}
b_0(\phi) & a_0^*(\phi) & 0 & 0 & 0 & \cdots \\ a_0(\phi) &
b_1(\phi) & a_1^*(\phi) & 0 & 0 & \cdots
\\ 0&a_1(\phi)&b_2(\phi)&a_2^*(\phi)&0&\cdots \\
\vdots&\vdots&\vdots&\vdots&\vdots&\ddots
\end{array}
\right)
\end{align*}
with the entries
\begin{align*}
a_n(\phi)&:\dom(a_n(\phi))\to\mathcal F_{n+1}(H),
\\ b_n(\phi)=(b_n(\phi))^*&:\dom(b_n(\phi))\to\mathcal F_n(H),\notag \\
a^*_n(\phi)=(a_n(\phi))^*&:\dom(a^*_n(\phi))\to\mathcal
F_n(H),\notag\\ \phi&\in H_+,~n\in\mathbb Z_+=0,1,\ldots.
\end{align*}
The inclusions $\dom(a_n(\phi))\subset\mathcal F_n(H)$,
$\dom(b_n(\phi))\subset\mathcal F_n(H)$, and
$\dom(a_n^*(\phi))\subset\mathcal F_{n+1}(H)$ hold for the
domains. We suppose these domains to contain $\mathcal F_n(H_+)$
and $\mathcal F_{n+1}(H_+)$, respectively.

Each matrix $\mathcal J(\phi)$ gives rise to a Hermitian operator
$J(\phi)$ in the space $\mathcal F(H)$: for
$\Phi=(\Phi_n)_{n=0}^\infty\in\dom(J(\phi))=\mathcal F_{fin}(H_+)$
we define \begin{align*}
(J(\phi)\Phi)_n&=a_{n-1}(\phi)\Phi_{n-1}+b_n(\phi)\Phi_n+a_n^*(\phi)\Phi_{n+1},\qquad
n\in\mathbb Z_+, \notag
\\ a_{-1}(\phi)&=0.
\end{align*}

Assume the following.
\begin{itemize}
\item[(a)] The operators $a_n(\phi)$
and $b_n(\phi)$, $\phi\in H_+$, $n\in\mathbb Z_+$, take real
spaces into real ones.
\item[(b)](smoothness) The restrictions $a_n(\phi)\upharpoonright\mathcal
F_n(H_+)$ and $b_n(\phi)\upharpoonright\mathcal F_n(H_+)$ act
continuously from $\mathcal F_n(H_+)$ to $\mathcal F_{n+1}(H_+)$
and $\mathcal F_n(H_+)$, respectively. The restrictions
$a_n^*(\phi)\upharpoonright\mathcal F_{n+1}(H_+)$ act continuously
from $\mathcal F_{n+1}(H_+)$ to $\mathcal F_n(H_+)$.
\item[(c)] The operators $J(\phi)$, $\phi\in H_+$, are essentially
selfadjoint and their closures $\tilde J(\phi)$, $\phi\in H_+$,
are strong commuting.
\item[(d)] The functions \begin{align*}
H_+\ni\phi &\mapsto a_n(\phi)\Phi_n\in\mathcal
F_{n+1}(H_+),~H_+\ni\phi\mapsto b_n(\phi)\Phi_n\in\mathcal
F_n(H_+),\\ H_+\ni\phi &\mapsto a^*_n(\phi)\Phi_{n+1}\in\mathcal
F_n(H_+),\qquad n\in\mathbb Z_+,
\end{align*}
are linear and continuous for all $\Phi_n\in\mathcal F_n(H_+)$,
$\Phi_{n+1}\in\mathcal F_{n+1}(H_+)$.
\item[(e)] (regularity) The real linear operators $V_n:\mathcal
F_n(H_+)\to\bigoplus_{j=0}^n\mathcal F_j(H_+)$ defined by the
equalities
\begin{align*}V_0=\Id_{\mathbb
C},~V_n(\phi_1\hat{\otimes}\cdots\hat\otimes\phi_n)&=
J(\phi_1)\ldots J(\phi_n)\Omega, \\ \phi_1,\ldots,\phi_n&\in
H_+,~n\in\mathbb N,\end{align*} are continuous. Furthermore, the
operators \begin{align*}\mathcal F_n(H_+)\ni F_n\mapsto
V_{n,n}F_n=(V_nF_n)_n\in\mathcal F_n(H_+),\qquad n\in\mathbb
N,\end{align*} are invertible.
\end{itemize}

We will call the family $J=(\tilde J(\phi))_{\phi\in H_+}$ of
operators a (commutative) Jacobi field if conditions (a)--(e) are
satisfied. Once again we should emphasize that the operators
$\tilde J(\phi)$ act in the Fock space $\mathcal F(H)$.

\subsection{Spectral theory of a Jacobi field}

It is possible to apply the projection spectral theorem, see~\cite{YuBYuK88} and~\cite{AP03}, to the
field $J=(\tilde J(\phi))_{\phi\in H_+}$. Here, we will only present the result of such an
application.

\begin{theorem}\label{thm_pst}
Given a Jacobi field $J$, there exist a Borel probability
measure $\rho$ on the space $H_-$ (the spectral measure) and a
vector-valued function $H_-\ni\xi\mapsto P(\xi)\in(\mathcal
F_{\mathrm{fin}}(H_+))'$ such that the following statements hold:
\begin{itemize}
\item[1.] For every $\xi\in H_-$, the vector
$P(\xi)=(P_n(\xi))_{n=0}^\infty\in(\mathcal F_{\mathrm{fin}}(H_+))'$, is a
generalized joint eigenvector of $J$ with eigenvalue $\xi$,
i.e.,
\begin{align}\label{gj_eig_def}\langle P(\xi),\tilde J(\phi)\Phi\rangle_H=\langle\xi,
\phi\rangle_H\langle P(\xi),\Phi\rangle_H,\qquad \phi\in
H_+,~\Phi\in\mathcal F_{\mathrm{fin}}(H_+).\end{align}
\item[2.]
After being extended by continuity to the whole of the space
$\mathcal F(H)$, the Fourier transform
\begin{align}\label{Fourier_J}\mathcal F(H)\supset\mathcal
F_{\mathrm{fin}}(H_+)\ni
\Phi&=(\Phi_n)_{n=0}^\infty\mapsto(I\Phi)(\xi)=\langle
\Phi,P(\xi)\rangle_H\notag
\\ &=\sum_{n=0}^\infty\langle
\Phi_n,P_n(\xi)\rangle_H\in L^2(H_-,d\rho)
\end{align}
becomes a unitary operator between $\mathcal F(H)$ and
$L^2(H_-,d\rho)$.
\item[3.]
The mapping $I$ takes every operator $\tilde J(\phi)$, $\phi\in
H_+$, to the operator of multiplication by the function
$H_-\ni\xi\mapsto\langle\xi,\phi\rangle_H\in\mathbb R$ in the
space $L^2(H_-,d\rho)$.
\end{itemize}
\end{theorem}
\begin{remark} The equality
\begin{align}\label{IV_m=monom} IV_nF_n=\langle\xi^{\otimes
n},F_n\rangle_H,\qquad F_n\in\mathcal F_n(H_+),~n\in\mathbb
Z_+,\end{align} holds true. Indeed, Assertion~3 of
Theorem~\ref{thm_pst} implies~(\ref{IV_m=monom}) for the vectors
\begin{align*}\sigma_n=\sum_{k=1}^l\lambda_k\phi_{1,k}\hat\otimes\cdots\hat\otimes
\phi_{n,k}&\in\mathcal F_n(H_+), \\ \lambda_k\in\mathbb
C,~\phi_{i,k}&\in H_+,~i=1,\ldots,n,~l\in\mathbb N.\end{align*} If
a sequence $(\sigma_n^i)_{i=0}^\infty$ of such vectors converges
to $F_n$ in the space $\mathcal F_n(H_+)$, then
\begin{align*}\langle\xi^{\otimes n},\sigma_n^i\rangle_H=IV_n\sigma_n^i\to
IV_nF_n\end{align*} in the space $L^2(H_-,d\rho)$. Since
$\langle\xi^{\otimes n},\sigma_n^i\rangle_H\to\langle\xi^{\otimes
n},F_n\rangle_H$ for each $\xi\in H_-$, the above formula implies
$IV_nF_n=\langle\xi^{\otimes n},F_n\rangle_H$.
\end{remark}

Now we have to recall some additional facts about the Fourier
transform~$I$.

Let $\mathcal P_n(H_-)$ denote the set of all continuous
polynomials on $H_-$ of degree $\le n$:  \begin{align*}H_-\ni \xi\mapsto
p_n(\xi)=\sum_{j=0}^n\langle\xi^{\otimes j},a_j\rangle_H\in
\mathbb C, \qquad a_j\in \mathcal F_j(H_+),~n\in \mathbb Z_+.\end{align*}

\begin{theorem}\label{IF_n(H)=calP_n}
The Fourier transform $I$ takes the set $\bigoplus_{j=0}^n\mathcal
F_j(H_+)\subset\mathcal F(H)$, $n\in \mathbb Z_+$, to the set
$\mathcal P_n(H_-)\subset L^2(H_-,d\rho)$ of continuous
polynomials on $H_-$ of degree $\le n$, i.e.,
\begin{align*}
I\left(\bigoplus_{j=0}^n\mathcal F_j(H_+)\right)=\mathcal
P_n(H_-),\qquad n\in \mathbb Z_+.
\end{align*}

The set $\mathcal P(H_-)=\bigcup_{n=0}^\infty\mathcal P_n(H_-)$ of
all continuous polynomials on $H_-$ is dense in $L^2(H_-,d\rho)$.
\end{theorem}

If $\dim H=\infty$, then $\bigoplus_{j=0}^n\mathcal F_j(H_+)$ is
not closed in $\mathcal F(H)$ and neither is $\mathcal P_n(H_-)$
closed in $L^2(H_-,d\rho)$. The closure of $\mathcal P_n(H_-)$ in
$L^2(H_-,d\rho)$ will be denoted by $\tilde{\mathcal P}_n(H_-)$.
The elements of $\tilde{\mathcal P}_n(H_-)$ are, by definition,
ordinary polynomials on $H_-$. Clearly,
\begin{align*}
I\left(\bigoplus_{j=0}^n\mathcal F_j(H)\right)=\tilde{\mathcal
P}_n(H_-),\qquad n\in \mathbb Z_+.
\end{align*}

The orthogonal decomposition $\mathcal
F(H)=\bigoplus_{n=0}^\infty\mathcal F_n(H)$ and the unitarity  of
$I$ imply the following orthogonal (chaotic) decomposition
of the space $L^2(H_-,d\rho)$:
\begin{align}\label{initial_dec}L^2(H_-,d\rho)&=\bigoplus_{n=0}^\infty
(L^2_n), \notag \\ (L^2_0)=\mathbb C,~(L^2_n)=I(\mathcal
F_n(H))&=\tilde{\mathcal P}_n(H_-)\ominus \tilde{\mathcal
P}_{n-1}(H_-),\qquad n\in\mathbb N.
\end{align}

\begin{remark}\label{prop_1.3}
Suppose $\mathcal H$ to be a nuclear space densely and
continuously embedded into $H$. In all the previous constructions,
it is possible to use the rigging
\begin{align*}\mathcal H'\supset H\supset\mathcal H\end{align*}
instead of the rigging~(\ref{rigging_H}). In this case, the family
$J$ consists of the operators $\tilde J(\phi)$, $\phi\in\mathcal
H$. The corresponding spectral measure $\rho$ is a Borel
probability measure on $\mathcal H'$.
\end{remark}

\subsection{Mapping of the spectral measure}

Consider a real separable Hilbert space $T_+$ and a rigging
\begin{align}\label{rigging T}T_-\supset T_0\supset T_+.\end{align}
As in the case of the rigging~(\ref{rigging_H}), the pairing
in~(\ref{rigging T}) can be extended to a pairing between
$\mathcal F_n(T_+)$ and $\mathcal F_n(T_-)$. The latter can, in
turn, be extended to a pairing between $\mathcal F_{fin}(T_+)$ and
$(\mathcal F_{fin}(T_+))'$. We use the notation
$\langle\cdot,\cdot\rangle_T$ for all of these pairings.

Let $K:T_+\to H_+$ be a linear continuous operator with
$\Ker(K)=\{0\}$ and suppose $\ran(K)$ to be dense in $H_+$. The
adjoint of $K$ with respect to~(\ref{rigging_H}) and~(\ref{rigging
T}) is a linear continuous operator $K^+:H_-\to T_-$ defined by
the equality
\begin{align*}\langle K^+\xi,f\rangle_T=\langle \xi,
Kf\rangle_H,\qquad \xi \in H_-,~f\in T_+.\end{align*}

\begin{lemma}
The kernel $\Ker(K^+)=\{0\}$. The range $\ran(K^+)$ is dense in
$T_-$.
\end{lemma}
\begin{proof}
Suppose $K^+\xi=0$ for some $\xi\in H_-$. This means $\langle
K^+\xi,f\rangle_T=\langle \xi,Kf\rangle_H=0$ for all $f\in T_+$.
Since $\ran(K)$ is dense in $H_+$, the latter implies $\xi=0$.
Thus $\Ker(K^+)=\{0\}$.

Next, we introduce a standard unitary $\mathbf I_T:T_-\to T_+$ by  the
formula
\begin{align*}(\mathbf I_T\omega,f)_{T_+}=\langle\omega,f\rangle_T,\qquad \omega\in T_-,~f\in
T_+.\end{align*} The equality
\begin{align*}(K^+\xi,\chi)_{T_-}=\langle K^+\xi,\mathbf I_T\chi\rangle_T
=\langle\xi,K\mathbf I_T\chi\rangle_H,\qquad \xi\in H_-,~\chi\in
T_-,\end{align*} holds true. If $(K^+\xi,\chi)_{T_-}=0$ for any
$\xi\in H_-$, then $K\mathbf I_T\chi=0$ and $\chi=0$. Thus
$\ran(K^+)$ is dense in $T_-$.
\end{proof}

Let $\mathcal B(H_-)$ stand for the Borel $\sigma$-algebra of the
space $H_-$. We denote by $\rho_K$ the image of the measure $\rho$
under the mapping $K^+$. By definition, $\rho_K$ is a probability
measure on the $\sigma$-algebra
\begin{align*}\mathcal C=\{\Delta\subset
T_-|(K^+)^{-1}(\Delta)\in\mathcal B(H_-)\},\end{align*}
($(K^+)^{-1}(\Delta)$ denoting the preimage of the set $\Delta$).
Clearly, the mapping $K^+$ is Borel-measurable, therefore
$\mathcal C$ contains the Borel $\sigma$-algebra of the space
$T_-$. If $K^+$ takes Borel subsets of $H_-$ to the Borel subsets
of $T_-$, then $\mathcal C$ coincides with the Borel
$\sigma$-algebra of $T_-$.

\section{Main results}\label{sec_JS}

Consider a Jacobi field $J=(\tilde J(\phi))_{\phi \in H_+}$ in the
Fock space $\mathcal F(H)$. The spectral measure $\rho$ of the
field $J$ is defined on $H_-$. The mapping $K^+$ takes $\rho$ to
the measure $\rho_K$ on $T_-$. The main objectives of this section
are:
\begin{itemize}
\item[1.] To obtain a family $J_K=(\tilde J_K(f))_{f\in T_+}$ of
commuting selfadjoint operators operators whose spectral measure
is  equal to $\rho_K$.
\item[2.] To obtain an analogue of the
decomposition~(\ref{initial_dec}) for the space
$L^2(T_-,d\rho_K)$.
\end{itemize}
We note that the family $J_K$ proves to satisfy
conditions~(a)--(d) of a Jacobi field. It is generally unclear
whether $J_K$ satisfies condition~(e).

The assumption $\Ker(K)=\{0\}$ is not essential. Indeed, the
measure $\rho_K$ proves to be lumped on the set of functionals
which equal zero on $\Ker(K)$. This set can be naturally
identified with $(\Ker(K)^\bot)'$. Thus we can always replace
$T_+$ with $\Ker(K)^\bot\subset T_+$.

\subsection{$\rho_K$ as the spectral measure}\label{rho_K spec meas}

Define the Hilbert space $T$ as the completion of $T_+$ with
respect to the scalar product
\begin{align*}(f_1,f_2)_T=(Kf_1,Kf_2)_H,\qquad f_1,f_2\in
T_+.\end{align*} The operator $K$ induces a unitary $\bar K:T\to
H$. We preserve the notations $K$ and $\bar K$ for the extensions
of $K$ and $\bar K$ to the complexified spaces $T_{+,c}$ and
$T_c$.

In the Fock space $\mathcal F(T)=\bigoplus_{n=0}^\infty \mathcal
F_n(T)$, consider the family $(\mathcal J_K(f))_{f\in T_+}$ of
operator-valued Jacobi matrices
\begin{align*}\mathcal J_K(f)=\left(
\begin{array}{cccccc}
\beta_0(f) & \alpha_0^*(f) & 0 & 0 & 0 & \cdots \\ \alpha_0(f) &
\beta_1(f) & \alpha_1^*(f) & 0 & 0 & \cdots
\\ 0&\alpha_1(f)&\beta_2(f)&\alpha_2^*(f)&0&\cdots \\
\vdots&\vdots&\vdots&\vdots&\vdots&\ddots
\end{array}
\right)\end{align*} with the entries
\begin{align*}
\alpha_n(f)=(\bar K^{\otimes(n+1)})^{-1}a_n(Kf)\bar K^{\otimes
n}&:\dom(\alpha_n(f))\to\mathcal F_{n+1}(T),\notag \\
\beta_n(f)=(\bar K^{\otimes n})^{-1}b_n(Kf)\bar K^{\otimes
n}&:\dom(\beta_n(f))\to\mathcal F_n(T),\notag \\
\alpha_n^*(f)=(\alpha_n(f))^*&:\dom(\alpha^*_n(f))\to\mathcal
F_n(T),\\ f&\in T_+,~n\in\mathbb Z_+,\end{align*} (recall that
$a_n(\phi)$ and $b_n(\phi)$ denote the entries of $\mathcal
J(\phi)$). The domains $\dom(\alpha_n(f))$, $\dom(\beta_n(f))$,
and $\dom(\alpha^*_n(f))$ contain $\mathcal F_n(T_+)$ and
$\mathcal F_{n+1}(T_+)$, respectively. As in the case of $\mathcal
J(\phi)$, each matrix $\mathcal J_K(f)$ gives rise to a Hermitian
operator $J_K(f)$ in the space $\mathcal F(T)$. The domain
$\dom(J_K(f))$ equals $\mathcal F_{fin}(T_+)$.

As we will further show, the operators $J_K(f)$, $f\in T_+$, are
essentially selfadjoint in the space $\mathcal F(T)$. Their
closures $\tilde J_K(f)$ are strongly commuting. Denote
$J_K=(\tilde J_K(f))_{f\in T_+}$.

\begin{theorem}\label{main_theorem}
Assume the restrictions \begin{align}\label{assump restr}a_n(Kf)
\upharpoonright(\ran(K))^{\hat\otimes
n},~b_n(Kf)&\upharpoonright(\ran(K))^{\hat\otimes n},\notag \\
a_n^*(Kf) &\upharpoonright(\ran(K))^{\hat\otimes(n+1)},\qquad f\in
T_+,~n\in\mathbb Z_+,
\end{align} to take values in
$(\ran(K))^{\hat\otimes(n+1)}$ and $(\ran(K))^{\hat\otimes n}$,
respectively. There exists  a vector-valued function
$T_-\ni\omega\mapsto Q(\omega)\in(\mathcal
F_{\mathrm{fin}}(T_+))'$ such that the following statements hold:
\begin{itemize}
\item[1.] For $\rho_K$-almost all $\omega\in T_-$, the vector
$Q(\omega)=(Q_n(\omega))_{n=0}^\infty\in(\mathcal
F_{\mathrm{fin}}(T_+))'$, is a generalized joint eigenvector of
the family $J_K$ with the eigenvalue $\omega$, i.e.,
\begin{align}\label{Q_is_eig}\langle Q(\omega),\tilde J_K(f)F\rangle_T=\langle\omega,f\rangle_T\langle
Q(\omega),F\rangle_T,\qquad F\in\mathcal
F_{\mathrm{fin}}(T_+).\end{align}
\item[2.]
After being extended by continuity to the whole of the space
$\mathcal F(T)$, the Fourier transform
\begin{align}\label{I_K_defin}
\mathcal F(T)\supset\mathcal F_{\mathrm{fin}}(T_+)\ni
F&=(F_n)_{n=0}^\infty\mapsto(I_KF)(\omega)=\langle
F,Q(\omega)\rangle_T\notag \\ &=\sum_{n=0}^\infty\langle
F_n,Q_n(\omega)\rangle_T\in L^2(T_-,d\rho_K)
\end{align}
becomes a unitary between $\mathcal F(T)$ and $L^2(T_-,d\rho_K)$.
\item[3.]
The mapping $I_K$ takes every operator $\tilde J_K(f)$, $f\in
T_+$, to the operator of multiplication by the function
$T_-\ni\omega\mapsto\langle\omega,f\rangle_T\in\mathbb R$ in the
space $L^2(T_-,d\rho_K)$.
\end{itemize}
\end{theorem}

\begin{proof}\textbf{Step 1.} First, we have to prove that
the operators $J_K(f)$, $f\in T_+$, are essentially selfadjoint
and their closures are strongly commuting. We define the operator
\begin{align}\label{calK_def}{\mathcal K}=\bigoplus_{n=0}^\infty
\bar K^{\otimes n}:\mathcal F(T)\to\mathcal F(H).\end{align} The
unitarity of $\bar K$ implies the unitarity of $\mathcal K$. A
straightforward calculation shows that
\begin{align}\label{J_K=KJ(Kf)K} J_K(f)={\mathcal K}^{-1}J(Kf){\mathcal
K},\qquad f\in T_+.\end{align} The operators $J(Kf)$ are
essentially selfadjoint and their closures are strongly commuting.
Since $\mathcal K$ is a unitary, the operators $J_K(f)$ possess
these properties, too.

\vspace{10pt}\textbf{Step 2.} Let us establish an isomorphism
between the spaces $L^2(T_-,d\rho_K)$ and $L^2(H_-,d\rho)$. For a
complex-valued function $G(\omega)$ on $T_-$, we define the
function
\begin{align}\label{U definit}(\mathcal UG)(\xi)=G(K^+\xi),\qquad
\xi\in H_-.\end{align} According to the definition of $\rho_K$,
the mapping $\mathcal U$ induces an isometric operator $U$ between
the spaces $L^2(T_-,d\rho_K)$ and $L^2(H_-,d\rho)$.

We have  $\ran(U)=L^2(H_-,d\rho)$. Indeed, consider an arbitrary
function $F(\xi)$ over $H_-$. Define
$G(\omega)=F((K^+)^{-1}\omega)$ if $\omega\in\ran(K^+)$ and
$G(\omega)=0$ otherwise. The equality \begin{align*}(\mathcal
UG)(\xi)=G(K^+\xi)=F(\xi),\qquad\xi\in H_-,\end{align*} holds
true. If $F\in L^2(H_-,d\rho)$, then $G\in L^2(T_-,d\rho_K)$. In
this case, the above equality yields $(UG)(\xi)=F(\xi)$.

As a result, we have the unitary $U:L^2(T_-,d\rho_K)\to
L^2(H_-,d\rho)$.

\vspace{10pt}\textbf{Step 3.} Consider the operator
\begin{align*}I_K=U^{-1}I\mathcal K:\mathcal F(T)\to
L^2(T_-,d\rho_K)\end{align*} with $I$ and $\mathcal K$ given
by~(\ref{Fourier_J}) and~(\ref{calK_def}), respectively. Since all
of its components are unitaries between the corresponding spaces,
$I_K$ is a unitary itself. Our next goal is to establish
representation~(\ref{I_K_defin}) for $I_K$.

Fix a vector $F\in\mathcal F_{\mathrm{fin}}(T_+)$. According
to~(\ref{Fourier_J}), the equality
\begin{align*}(I\mathcal KF)(\xi)=\langle\mathcal
KF,P(\xi)\rangle_H=\sum_{n=0}^\infty\langle F,(K^+)^{\otimes
n}P_n(\xi)\rangle_H,\qquad \xi\in H_-,\end{align*} holds true.
Define $Q(K^+\xi)=((K^+)^{\otimes
n}P_n(\xi))_{n=0}^\infty\in(\mathcal F_{\mathrm{fin}}(T_+))'$.
Note that $Q(K^+\xi)$ is well-defined because $K^+$ is
monomorphic. Evidently,
\begin{align}\label{<KF,P>=<F,Q>}(I\mathcal KF)(\xi)=\langle
F,Q(K^+\xi)\rangle_T,\qquad \xi\in H_-.\end{align} The application
of $U^{-1}$ to~(\ref{<KF,P>=<F,Q>}) yields
representation~(\ref{I_K_defin}) for the unitary
$I_K=U^{-1}I\mathcal K$.

The proof of Theorem~\ref{main_theorem} will be complete if we
show that Statements~1 and~3 hold for the function $Q(\omega)$.

\vspace{10pt}\textbf{Step 4.} Let us prove Statement~1. As before,
we fix a vector $F\in\mathcal F_{\mathrm{fin}}(T_+)$. Due to the
assumption on the restrictions~(\ref{assump restr}), the vector
$\tilde J_K(f)F$ belongs to $\mathcal F_{\mathrm{fin}}(T_+)$.
Hence the right-hand side of~(\ref{Q_is_eig}) is well-defined.

Formulas~(\ref{gj_eig_def}),~(\ref{J_K=KJ(Kf)K}),
and~(\ref{<KF,P>=<F,Q>}) imply
\begin{align*}
\langle Q(\omega),\tilde J_K(f)F\rangle_T&=\overline{\langle
\tilde J_K(f)F,Q(\omega)\rangle}_T \\ &=\overline{(I_K\tilde
J_K(f)F)}(\omega)\\&=\overline{(U^{-1}I\mathcal K\tilde
J_K(f)F)}(\omega) \\ &=\overline{(U^{-1}I\mathcal K\mathcal
K^{-1}\tilde J(Kf)\mathcal KF)}(\omega)\\ &
=\overline{(U^{-1}I\tilde J(Kf)\mathcal KF)}(\omega) \\
&=(U^{-1}\overline{\langle \tilde
J(Kf)\mathcal KF,P(\cdot)\rangle}_H)(\omega) \\
&=(U^{-1}(\overline{\langle Kf,\cdot\rangle_H\langle\mathcal
KF,P(\cdot)\rangle}_H))(\omega)
\\ &=(U^{-1}\overline{\langle
Kf,\cdot\rangle}_H)(\omega)(U^{-1}\overline{\langle \mathcal
KF,P(\cdot)\rangle}_H)(\omega) \\ &=\langle
f,\omega\rangle_T\overline{(U^{-1}I\mathcal KF)}(\omega)
\\ &=\langle f,\omega\rangle_T\overline{(I_KF)}(\omega)\\&=\langle
f,\omega\rangle_T\langle Q(\omega),F\rangle_T
\end{align*}
for $\rho_K$-almost all $\omega\in T_-$ (overbars denote complex
conjugacy). This proves Statement~1.

Statement~3 is a direct consequence of~(\ref{I_K_defin})
and~(\ref{Q_is_eig}).
\end{proof}

\begin{remark}
While proving the theorem, we showed that the mapping~(\ref{U
definit}) induces a unitary $U:L^2(T_-,d\rho_K)\to
L^2(H_-,d\rho)$. We also obtained an explicit formula for the
Fourier transform $I_K$. Namely,
\begin{align}\label{I_K expl} I_K=U^{-1}I\mathcal K\end{align} with $I$ and $\mathcal
K$ given by~(\ref{Fourier_J}) and~(\ref{calK_def}), respectively.
\end{remark}
\begin{remark}
As mentioned above, it is generally unclear whether $J_K$
satisfies condition (e) in the definition of a Jacobi field.
However, if the operator $K$ is invertible, then $J_K$ does
satisfy (e) and hence is a Jacobi field. This field is isomorphic
to the initial field $J$.
\end{remark}

\subsection{Orthogonal (chaotic) decomposition of the space\\ $L^2(T_-,d\rho_K)$}

This subsection aims to obtain an analogue of the
decomposition~(\ref{initial_dec}) for the space
$L^2(T_-,d\rho_K)$. If $J_K$ proves to be a Jacobi field, then
Theorem~\ref{IF_n(H)=calP_n} is applicable. Otherwise, an analogue
of~(\ref{initial_dec}) for $L^2(T_-,d\rho_K)$ may be obtained with
the help  of Theorem~\ref{UQn=Pn}.

Further considerations do not require any assumptions on the
restrictions~(\ref{assump restr}). Theorem~\ref{UQn=Pn} below is
applicable  to a Jacobi field which does not satisfy the
assumption of Theorem~\ref{main_theorem}. In this case, the
unitary $I_K$ should be defined by formula~(\ref{I_K expl}).

Let $\mathcal Q_n(T_-)$ denote the set of all continuous
polynomials
\begin{align}\label{pol_on_T}T_-\ni \omega\mapsto
q_n(\omega)=\sum_{j=0}^n\langle\omega^{\otimes j},c_j\rangle_T\in
\mathbb C,\qquad c_j\in \mathcal F_j(T_+),~n\in \mathbb
Z_+,\end{align} on  $T_-$ of degree $\le n$. As will be shown
below, the inclusion
\begin{align}\label{Qn subset L2}\mathcal Q_n(T_-)\subset
L^2(T_-,d\rho_K)\end{align} holds. The closure of $\mathcal
Q_n(T_-)$ in $L^2(T_-,d\rho_K)$ will be denoted by
$\tilde{\mathcal Q}_n(T_-)$. The elements of $\tilde{\mathcal
Q}_n(T_-)$ are ordinary polynomials on  $T_-$.

\begin{theorem}\label{UQn=Pn}
The unitary $I_K$ takes the set $\bigoplus_{j=0}^n\mathcal
F_j(T)\subset\mathcal F(T)$, $n\in \mathbb Z_+$, to the set
$\tilde{\mathcal Q}_n(T_-)\subset L^2(T_-,d\rho_K)$ of ordinary
polynomials on $T_-$, i.e.,
\begin{align}\label{I_KsumFn=Qn}
I_K\left(\bigoplus_{j=0}^n\mathcal F_j(T)\right)=\tilde{\mathcal
Q}_n(T_-),\qquad n\in \mathbb Z_+.
\end{align}

The set $\mathcal Q(T_-)=\bigcup_{n=0}^\infty\mathcal Q_n(T_-)$ of
all continuous polynomials on $T_-$ is dense in
$L^2(T_-,d\rho_K)$.
\end{theorem}
\begin{proof}
\textbf{Step 1.} First, we have to prove inclusion~(\ref{Qn subset
L2}). The application of the mapping~(\ref{U definit}) to the
polynomial~(\ref{pol_on_T}) yields
\begin{align*}(\mathcal Uq_n)(\xi)=q_n(K^+\xi)&=\sum_{j=0}^n \langle (K^+\xi)^{\otimes
j},c_j\rangle_T\notag \\ &=\sum_{j=0}^n \langle (K^+)^{\otimes
j}\xi^{\otimes j},c_j\rangle_T=\sum_{j=0}^n \langle \xi^{\otimes
j},K^{\otimes j}c_j\rangle_H.
\end{align*}
The expression in the right hand side of this formula is a
continuous polynomial with the coefficients $a_j=K^{\otimes
j}c_j\in \mathcal F_j(H_+)$. According to
Theorem~\ref{IF_n(H)=calP_n}, this polynomial belongs to the space
$L^2(H_-,d\rho)$. Therefore $q_n(\omega)$ belongs to the space
$L^2(T_-,d\rho_K)$. The latter proves inclusion~(\ref{Qn subset
L2}).

\vspace{10pt}\textbf{Step 2.} Let us prove
equality~(\ref{I_KsumFn=Qn}). Formula~(\ref{I_K expl}) and
Theorem~\ref{IF_n(H)=calP_n} yield
\begin{align*}I_K\left(\bigoplus_{j=0}^n\mathcal
F_j(T)\right)&=U^{-1}I\mathcal K\left(\bigoplus_{j=0}^n\mathcal
F_j(T)\right) \\ &=U^{-1}I\left(\bigoplus_{j=0}^n\mathcal
F_j(H)\right)=U^{-1}\tilde{\mathcal P}_n(H_-),\qquad n\in \mathbb
Z_+.\end{align*} The proof of equality~(\ref{I_KsumFn=Qn}) will be
complete if we show that $\tilde{\mathcal
P}_n(H_-)=U\tilde{\mathcal Q}_n(T_-)$, $n\in\mathbb Z_+$.

As explained above, each function $Uq_n(\xi)=q_n(K^+\xi)$,
$q_n(\omega)\in\mathcal Q_n(T_-)$, belongs to $\mathcal P_n(H_-)$.
Thus it is only necessary to prove that such functions are dense
in $\tilde{\mathcal P}_n(H_-)$.

\vspace{10pt}\textbf{Step 3.} It suffices to approximate a
monomial $\langle\xi^{\otimes m},a_m \rangle_H$, $a_m\in\mathcal
F_m(H_+)$, $m=1,\ldots,n$, with the elements of $U\mathcal
Q_n(T_-)$.

Fix $\epsilon>0$. Since $\ran(K)$ is dense in $H_+$, there exists
a vector
\begin{align*}s_{m,\epsilon}=\sum_{k=1}^l\lambda_kf_{1,k}\hat\otimes\cdots\hat\otimes
f_{m,k}&\in\mathcal F_m(T_+), \\ \lambda_k\in\mathbb
C,~f_{i,k}&\in T_+,~i=1,\ldots,m,~l\in\mathbb N,\end{align*} such
that
\begin{align*}\|a_m-K^{\otimes n}s_{m,\epsilon}\|_{\mathcal
F(H_+)}=\left\|a_m-\sum_{k=1}^l\lambda_k
Kf_{1,k}\hat\otimes\cdots\hat\otimes Kf_{m,k}\right\|_{\mathcal
F(H_+)}<\epsilon.\end{align*} Taking equality~(\ref{IV_m=monom})
into account, we conclude that the monomial
$\langle\omega^{\otimes m},s_{m,\epsilon}\rangle_T\in\mathcal
Q_n(T_-)$ satisfies the estimate
\begin{align*}\|\langle\xi^{\otimes
m},a_m\rangle_H&-(U\langle\cdot^{\otimes
m},s_{m,\epsilon}\rangle_T)(\xi)\|_{L^2(H_-,d\rho(\xi))} \\
&=\|\langle\xi^{\otimes m},a_m\rangle_H-\langle (K^+)^{\otimes
m}\xi^{\otimes m},s_{m,\epsilon} \rangle_T\|_{L^2(H_-,d\rho(\xi))}
\\ &=\|\langle\xi^{\otimes m},a_m\rangle_H-\langle \xi^{\otimes
m},K^{\otimes m}s_{m,\epsilon}\rangle_H\|_{L^2(H_-,d\rho(\xi))} \\
&=\|I^{-1}(\langle\xi^{\otimes m},a_m\rangle_H-\langle
\xi^{\otimes m},K^{\otimes m}s_{m,\epsilon}\rangle_H)\|_{\mathcal
F(H)} \\ &=\|I^{-1}(IV_ma_m-IV_mK^{\otimes
m}s_{m,\epsilon})\|_{\mathcal F(H)}
\\ &=\|V_m(a_m-K^{\otimes
m}s_{m,\epsilon})\|_{\mathcal F(H)}
\\ &\le\|V_m(a_m-K^{\otimes
m}s_{m,\epsilon})\|_{\mathcal F(H_+)} \\
&\le\|V_m\|\,\|a_m-K^{\otimes m}s_{m,\epsilon}\|_{\mathcal
F(H_+)}<\|V_m\|\epsilon.
\end{align*}
Thus we have approximated $\langle\xi^{\otimes m},a_m \rangle_H$
with the functions $(U\langle\cdot^{\otimes
m},s_{m,\epsilon}\rangle_T)(\xi)\in U\mathcal Q_n(T_-)$.

\vspace{10pt}\textbf{Step 4.} Let us prove the last assertion of
Theorem~\ref{UQn=Pn}. Due to the unitarity of $I_K$,
\begin{align*}
L^2(T_-,d\rho_K)=I_K(\mathcal F(T))&=(I_K(\mathcal
F_{fin}(T)))^{\sim}=\left(\bigcup_{n=0}^\infty
I_K\left(\bigoplus_{m=0}^n\mathcal
F_m(T)\right)\right)^{\sim}\\&=\left(\bigcup_{n=0}^\infty\tilde{\mathcal
Q}_n(T_-)\right)^{\sim}=\left(\bigcup_{n=0}^\infty\mathcal
Q_n(T_-)\right)^{\sim},
\end{align*}
(tilde stands for the closure in the corresponding space). Thus
$\mathcal Q(T_-)$ is dense in $L^2(T_-,d\rho_K)$.
\end{proof}

We can now construct the~(\ref{initial_dec})-type decomposition
for the space $L^2(T_-,d\rho_K)$:
\begin{align*}L^2(T_-,d\rho_K)&=\bigoplus_{n=0}^\infty
(L^2_n)_K,\notag \\ (L^2_0)_K=\mathbb C,~(L^2_n)_K=I_K(\mathcal
F_n(T))&=\tilde{\mathcal Q}_n(T_-)\ominus \tilde{\mathcal
Q}_{n-1}(T_-),\qquad n\in\mathbb N.
\end{align*}

\section{Examples}

Let us make some remarks concerning the space $T$ and the Fourier
transform of the measure $\rho_K$,
\begin{align*}\hat\rho_K(f)=\displaystyle\int_{T_-}e^{i\langle\omega,f\rangle_T}d\rho_K(\omega),
\qquad f\in T_+.\end{align*}  We will be using these remarks in
our further considerations.

\begin{remark}
Since $K:T_+\to H_+$ is continuous and since the embedding
$H_+\hookrightarrow H$ is continuous, we easily conclude that
$T_+$ is continuously embedded into $T$. Furthermore, $T_+$ is a
dense subset of $T$. Thus we can use $T$ as the zero space in the
chain~(\ref{rigging T}), i.e., we can assume $T_0=T$.
\end{remark}

\begin{remark}\label{K^+=K^-1}
The set $\ran(K)\subset H_+\subset H_-$ is dense in $H_-$.
Assuming $T_0=T$, one can prove that the restriction
$K^+\upharpoonright\ran(K):\ran(K)\to T_-$ coincides with the
mapping $K^{-1}:\ran(K)\to T_+\subset T_-$.
\end{remark}

\begin{remark}\label{Fourier_measures} Consider the Fourier transform
\begin{align*}\hat\rho(\phi)=\int_{H_-}e^{i\langle\xi,\phi\rangle_H}d\rho(\xi),
\qquad \phi\in H_+,\end{align*} of the measure $\rho$. By the
definition of $\rho_K$, we have:
\begin{align*}\hat\rho_K(f)=\int_{T_-}e^{i\langle\omega,f\rangle_T}d\rho_K(\omega)
&=\int_{H_-}e^{i\langle K^+\xi,f\rangle_H}d\rho(\xi) \\
&=\int_{H_-}e^{i\langle\xi,Kf\rangle_H}d\rho(\xi)
=\hat\rho(Kf),\qquad f\in T_+.\end{align*} Thus the Fourier
transform $\hat\rho_K(f)$ of the measure $\rho_K$ satisfies the
equality
\begin{align*}\hat\rho_K(f)=\hat\rho(Kf),\qquad f\in T_+.\end{align*}
\end{remark}

We will now apply the results of Section~\ref{sec_JS} to some
classical Jacobi fields.

\begin{example}\label{ex Gauss}
Suppose $J$ to be the classical free field, see
e.g.~\cite{YuBYuK88},
\cite{YuBVLEL95},~\cite{EL95},~\cite{YuB98IEO},
and~\cite{YuB98MF}. In this case,
\begin{align*}a_n(\phi)\Phi_n=(\sqrt{n+1}\phi)\hat\otimes \Phi_n,~b_n(\phi)\Phi_n=0,\qquad
\Phi_n\in\mathcal F_n(H),~\phi\in H_+,~n\in\mathbb
Z_+.\end{align*} Clearly, the assumption of
Theorem~\ref{main_theorem} on the restrictions \eqref{assump
restr} is now automatically satisfied for any operator $K$ under
consideration.

The spectral measure $\rho$ of the field $J$ is the standard
Gaussian measure $\gamma$ on $H_-$. Its Fourier transform is given
by the formula \begin{align*}\hat
\rho(\phi)=\hat\gamma(\phi)=\exp\left(-\frac12\,\|\phi\|_H^2\right),\qquad
\phi\in H_+.\end{align*} According to
Remark~\ref{Fourier_measures}, the Fourier transform of $\rho_K$
is given by the formula
\begin{align*}\hat\rho_K(f)=\hat\rho(Kf)=\exp\left(-\frac12\,\|Kf\|_H^2\right)=\exp\left(-\frac12\,\langle
K^+Kf,f\rangle_T\right),\qquad f\in T_+,\end{align*} (since $H_+$
is a subset of $H_-$, the operator $K^+K:T_+\to T_-$ is
well-defined). This means $\rho_K$ is the Gaussian measure on
$T_-$ with the correlation operator $K^+K$. Notice that, in the
case where $T_0=T$, the Fourier transform of $\rho_K$ may be
written down in the form
\begin{align*}\hat\rho_K(f)=\exp\left(-\frac12\,\|f\|_T^2\right),\qquad f\in T_+,\end{align*} i.e.,
$\rho_K$ is the standard Gaussian measure on $T_-$.

Applying Theorem~\ref{main_theorem} to the classical free field
$J$, we obtain the family $J_K$ whose spectral measure is
$\rho_K$.
\end{example}

In what follows, we assume $T_0=T$.

\begin{example}
Let $H$ be $L^2(\mathbb R,dx)$ and let $H_+$ and $T_+$ be the
Sobolev spaces\linebreak  $W_2^1(\mathbb R,(1+x^2)\,dx)$  and $W_2^1(\mathbb
R,dx)$, respectively. Suppose $J$ to be the Poisson field, see
e.g.~\cite{EL95},~\cite{YuB98IEO},~\cite{YuB98MF},
and~\cite{YuB00IDA}. In this case,
\begin{align*}a_0(\phi)\Phi_0&=\Phi_0\phi,~b_0(\phi)\Phi_0=0,
\\a_n(\phi)\Phi_n&=(\sqrt{n+1}\phi)\hat\otimes\Phi_n, \\ b_n(\phi)\Phi_n&=
(b(\phi)\otimes\Id_{H_+}\otimes\cdots\otimes\Id_{H_+})\Phi_n \\
&+(\Id_{H_+}\otimes
b(\phi)\otimes\Id_{H_+}\otimes\cdots\otimes\Id_{H_+})\Phi_n+\cdots
\\ &+(\Id_{H_+}\otimes \cdots\otimes\Id_{H_+}\otimes
b(\phi))\Phi_n,
\\ &\Phi_0\in\mathcal F_0(H_+),~\Phi_n\in\mathcal F_n(H_+),~\phi\in
H_+,~n\in\mathbb N.\end{align*} Here, $b(\phi)$ the operator of
multiplication by the function $\phi(x)$ in the space $H$.

The space $H_-$ coincides with the negative Sobolev space
$W_2^{-1}(\mathbb R,(1+x^2)\,dx)$. The spectral measure $\rho$ of
the field $J$ is equal to the centered Poisson measure $\pi$ with
the intensity $dx$. The Fourier transform of $\rho$ is given by
the formula
\begin{align*}\hat\rho(\phi)=\hat\pi(\phi)
=\exp\left(\int_\mathbb
R(e^{i\phi(x)}-1-i\phi(x))\,dx)\right),\qquad \phi\in
H_+.\end{align*}

Suppose $K:T_+\to H_+$ to be the operator of multiplication by the
function $\kappa(x)=e^{-x^2}$. One can easily verify that $K$ is
bounded and $\Ker(K)=\{0\}$. The range $\ran(K)$ is dense in $H_+$
because it contains all the smooth compactly supported  functions.
On the other hand,  $\ran(K)\ne H_+$ because e.g. the function
$\psi(x)=(1+x^2)^{-2}\in H_+$ does not belong to $\ran(K)$.
Clearly, the field $J$ and the operator $K$ satisfy the
assumptions of Theorem~\ref{main_theorem}.

The space $T_-$ is the dual of the space $W_2^1 (\mathbb R,dx)$
with respect to the zero space $T=L^2(\mathbb R,e^{-2x^2}\,dx)$.
Evidently, one may realize $T_-$ as the dual space of $W_2^1
(\mathbb R,dx)$ with respect to the zero space $L^2(\mathbb
R,dx)$, in which case $T_-$ is the usual negative Sobolev space
$W_2^{-1}(\mathbb R,dx)$.

According to Remark~\ref{K^+=K^-1}, the operator $K^+:H_-\to T_-$
is equal to the extension by continuity of the mapping
\begin{align*}H_-\supset\ran(K)\ni\xi(x)\mapsto e^{x^2}\xi(x)\in
T_-.\end{align*} According to Remark~\ref{Fourier_measures}, the
Fourier transform of $\rho_K$ is given by the formula
\begin{align*}\hat\rho_K(f(x))&=\hat\rho(e^{-x^2}f(x)) \\ &=
\exp\left(\int_\mathbb R\left(e^{ie^{-x^2}f(x)}-1-ie^{-x^2}
f(x)\right)\,dx\right),\qquad f\in T_+.\end{align*} Applying
Theorem~\ref{main_theorem}, we obtain the family $J_K$ whose
spectral measure is  $\rho_K$.
\end{example}

\begin{example}
As before, let $H$ be  $L^2(\mathbb R,dx)$. Let $H_+$ and $T_+$
equal $W_2^1\bigl(\mathbb R,e^{\frac{x^2}2}\,dx\bigr)$ and
$W_2^2\bigl(\mathbb R,e^{\frac{x^2}2}\,dx\bigr)$, respectively.
Suppose $J$ to be the Poisson field.

Define the operator $K:T_+\to H_+$ as the extension by continuity
of the mapping
\begin{align*}\mathrm C_0^\infty(\mathbb R)\ni p(x)\mapsto e^{-\frac{x^2}2}\frac{dp(x)}{dx}\in
H_+\end{align*} ($\mathrm C_0^\infty(\mathbb R)$ stands for the
set of all smooth compactly  supported functions on $\mathbb R$).
Evidently, $K$ is bounded and $\Ker(K)=\{0\}$.

\begin{lemma}
The range $\ran(K)$ is dense in $H_+$.
\end{lemma}
\begin{proof}
Fix $q\in H_+$ and assume $(Kp(x),q(x))_{H_+}=0$ for an arbitrary
$p\in\mathrm C_0^\infty(\mathbb R)$. Our goal is to show that
$q=0$.

The equality
\begin{align*}(Kp,q)_{H_+}=\int_\mathbb
R\left(-\frac{d^3p(x)}{dx^3}+x\frac{d^2p(x)}{dx^2}+2\frac{dp(x)}{dx}\right)q(x)\,dx\end{align*}
holds. Consider the differential expression
\begin{align*}\mathcal L=-\frac{d^3}{dx^3}+x\frac{d^2}{dx^2}+2\frac{d}{dx}.\end{align*}
 Let $\mathcal L^+$ denote the
adjoint expression. Since
\begin{align*}(Kp,q)_{H_+}=(\mathcal Lp,q)_H=0\end{align*}
for an arbitrary $p\in\mathrm C_0^\infty(\mathbb R)$, the function
$q$ is a generalized solution of the differential equation
$\mathcal L^+y=0$. Calculating $\mathcal L^+$ and applying
Theorem~6.1 from Chapter~16 of~\cite{YuBGUZSh90}, we conclude that
$q$ is indeed a classical solution of the equation
\begin{align*}\frac{d^3y(x)}{dx^3}+x\frac{d^2y(x)}{dx^2}=0.\end{align*}

The general solution of the above equation is
\begin{align*}y(x)=c_1\int_0^x\int_0^te^{-\frac{s^2}2}ds\,dt+c_2x+c_3,
\qquad c_1,c_2,c_3\in\mathbb C.
\end{align*}
Assume $y=q\in H_+$. In this case, the limits
\begin{align*}
\lim_{x\to\infty}\frac{dy(x)}{dx}&=c_1\lim_{x\to\infty}\int_0^xe^{-\frac{t^2}2}dt+c_2=c_1\sqrt\frac\pi2+c_2,
\\
\lim_{x\to-\infty}\frac{dy(x)}{dx}&=-c_1\lim_{x\to-\infty}\int_x^0e^{-\frac{t^2}2}dt+c_2=-c_1\sqrt\frac\pi2+c_2\end{align*}
must equal $0$. Evidently, the latter implies $c_1=c_2=c_3=0$.
Thus $q=0$.
\end{proof}

The Poisson field $J$ and the operator $K$ do not satisfy the
assumptions of Theorem~\ref{main_theorem}. However,
Theorem~\ref{UQn=Pn} is applicable now.

The space $H_-$ is  the negative Sobolev space
$W_2^{-1}\bigl(\mathbb R,e^\frac{x^2}2\,dx\bigr)$, while $T_-$ may
be realized as the dual of the space $W_2^2 (\mathbb
R,e^\frac{x^2}2\,dx)$ with respect to the zero space $L^2(\mathbb
R,dx)$. In this case, $T_-$ is the usual negative Sobolev space
$W_2^{-2}(\mathbb R,e^\frac{x^2}2\,dx)$. According to
Remark~\ref{Fourier_measures}, the Fourier transform of $\rho_K$
is given by the formula
\begin{align*}\hat\rho_K(f(x))&
=\hat\rho\left(e^{-\frac{x^2}2}\,\frac{df(x)}{dx}\right) \\ &=
\exp\left(\int_{\mathbb R}\left(\exp\left(ie^{-\frac{x^2}2}\,\frac
{df(x)}{dx}\right)-1-ie^{-\frac{x^2}2}\frac
{df(x)}{dx}\right)\,dx\right),\qquad f\in T_+.\end{align*}
Applying Theorem~\ref{UQn=Pn}, we obtain a
(\ref{initial_dec})-type decomposition for the space \linebreak
$L^2(T_-,d\rho_K)$.
\end{example}

In a forthcoming paper, we are going to discuss in detail the case
of the fractional Brownian motion, which is an important example
of a Gaussian measure with a non-trivial correlation operator.

\noindent Yurij M. Berezansky, Institute of  Mathematics,  National Academy of Sciences of Ukraine,  3
Tereshchenkivs'ka,  01601 Kiev,  Ukraine\\
\texttt{berezan@mathber.carrier.kiev.ua}\\[2mm]
Eugene W.  Lytvynov, Department of Mathematics,
University of Wales Swansea, Singleton Park, Swansea SA2 8PP, U.K. \\
 \texttt{e.lytvynov@swansea.ac.uk}\\[2mm]
Artem~D.~Pulemyotov, Department of Mathematics and Mechanics,  Kiev
National T. Shevchenko University,  64 Volodymyrs'ka,  01033
Kiev,  Ukraine\\
\texttt{pulen@i.kiev.ua}

\end{document}